\documentclass{amsart}
\usepackage{amssymb,amsmath,latexsym}
\newtheorem{theorem}{Theorem}
\newcommand{\bt}{\begin{theorem}}
\newcommand{\et}{\end{theorem}}
\newtheorem{lemma}{Lemma}
\newcommand{\bl}{\begin{lemma}}
\newcommand{\el}{\end{lemma}}
\newcommand{\bpf}{\begin{proof}}
\newcommand{\epf}{\end{proof}}
\newcommand{\beq}{\begin{equation}}
\newcommand{\eeq}{\end{equation}}
\newcommand{\benum}{\begin{enumerate}}
\newcommand{\eenum}{\end{enumerate}}
\newcommand{\N}{\ensuremath{ \mathbf N }}
\newcommand{\Z}{\ensuremath{\mathbf Z}}
\newcommand{\card}{\text{card}}

\title[Bases for binary linear forms]{Representation functions of bases for binary linear forms}

\date{\today}



\author{Melvyn B. Nathanson}
\thanks{This work was supported
in part by grants from the NSA Mathematical Sciences Program and
the PSC-CUNY Research Award Program}
\address{Department of Mathematics, Lehman College (CUNY), 
Bronx, NY 10468,
and School of Mathematics, Institute for Advanced Study,
Princeton, NJ 08540}
\email{melvyn.nathanson@lehman.cuny.edu, melvyn@ias.edu}

\dedicatory{To Jean-Marc Deshouillers on his 60th birthday}

\begin{document}

\begin{abstract}
Let $F(x_1,\ldots,x_m) = u_1 x_1 + \cdots + u_mx_m$ be a linear form with nonzero, 
relatively prime integer coefficients $u_1, \ldots, u_m$.  For any set $A$ of integers, 
let $F(A) = \{F(a_1,\ldots,a_m) : a_i \in A \text{ for $i=1,\ldots,m$}\}.$  
The \emph{representation function} associated with the form $F$ is 
\[
R_{A,F}(n) = \card\left( \{ (a_1,\ldots,a_m)\in A^m:  F(a_1,\ldots, a_m) = n \} \right).
\]
The set $A$ is a \emph{basis with respect to $F$ for almost all integers} if the set $\Z \setminus F(A)$ has asymptotic density zero.  Equivalently, the representation function of a basis for almost all integers is a function $f:\Z \rightarrow \N_0\cup\{\infty\}$ such that $f^{-1}(0)$ has density zero.  Given such a function, the inverse problem for bases is to construct a set $A$ whose representation function is $f$.  In this paper the inverse problem is solved for binary linear forms.
\end{abstract}

\maketitle

\section{Direct and inverse problems for representation functions}

Let \N, $\N_0,$ and \Z\ denote the sets of positive integers, nonnegative integers, and integers, respectively.

In additive number theory, a classical direct problem is to describe the integers that can be represented as the sum of a bounded number of elements of a fixed set $A$ of integers.  For example, describe the integers that are sums of two primes or three squares or four cubes.  Given a set $A$ and a positive integer $m,$ we associate various representation functions to the \emph{sumset} $mA = \{a_1 + \cdots + a_m : a_i \in A \text{ for $i=1,\ldots, m$} \}.$  The two most important are the \emph{ordered representation function}
\[
R_{A,m}(n) = \{ (a_1,\ldots, a_m) \in A^m : \sum_{i=1}^m a_i = n \}
\]
and the \emph{unordered representation function}
\[
r_{A,m}(n) = \{ (a_1,\ldots, a_m) \in A^m : \sum_{i=1}^m a_i = n     \text{ and }  a_1 \leq \cdots \leq a_m \}.
\]
Representation functions\footnote{It is semantically ironic that the unordered representation function counts only linearly ordered representations while the ordered representation function counts unordered representations.}
are functions from \Z\ into the set $\N_0 \cup \{\infty\}.$  

Let $A$ be a set of integers.  The \emph{counting function} $A(x_1,x_2)$ counts the number of elements of $A$ between $x_1$ and $x_2,$ that is,
\[
A(x_1,x_2) = \sum_{\substack{a \in A \\ x_1 \leq a \leq x_2} } 1.
\]
The set $A$ has \emph{upper asymptotic density} $d_U(A) = \alpha$ if $\limsup_{x\rightarrow\infty} A(-x,x)/(2x+1) = \alpha.$  
The set $A$ has \emph{asymptotic density} $d(A) = \alpha$ if $\lim_{x\rightarrow\infty} A(-x,x)/(2x+1) = \alpha.$ 
If $S$ is an infinite set of integers and $W$ is a subset of $S$, then the set $W$ has \emph{relative asymptotic density} $d(W,S) = \alpha$ if $\lim_{x\rightarrow\infty} W(-x,x)/S(-x,x) = \alpha.$ 

The set $A$ is called a \emph{basis of order $m$ for $S$} if $S \subseteq mA,$ that is, if $R_{A,m}(n) \geq 1$ for every element of $S$, and a \emph{basis of order $m$ for almost all $S$} if $\{ n \in S : R_{A,m}(n) = 0\}$ has relative asymptotic density zero.

Nathanson~\cite{nath78d} introduced a class of inverse problems for the representation functions of sets of integers and nonnegative integers.   Let $m \geq 2.$  If $A$ and $B$ are sets of integers such that $R_{A,m}(n) = R_{B,m}(n)$ for all $n\in \Z,$ then does $A=B$?  He proved that the answer is  ``yes'' if $A$ and $B$ are sets of nonnegative integers, but if  $R_{A,m}(n) = R_{B,m}(n)$  only for all sufficiently large integers $n,$ then $A$ and $B$ are not necessarily equal, but the structures of $A$ and $B$ can be described explicitly.

These results suggest a second class of inverse problems for representation functions.  Given any function $f: \Z \rightarrow \N_0 \cup \{\infty\}$ and an integer $m\geq 2$, does there exist a set $A$ such that $R_{A,m}(n) = f(n)$ for all $n\in \Z$ or $r_{A,m}(n) = f(n)$ for all $n\in \Z$?  Describe all sets $A$ such that $R_{A,m} = f$ or $r_{A,m} = f.$  
Nathanson~\cite{nath04a,nath05a}  proved that \emph{every}  function $f: \Z \rightarrow \N_0 \cup \{\infty\}$ such that  $f^{-1}(0)$ is finite is the unordered representation function for infinitely many sets of integers.  In particular, there exist \emph{unique representation bases} for the integers, that is, sets $A \subseteq \Z$ such that $r_{A,m}(n) = 1$ for all $n\in \Z$ (For $m=2,$ see {\L}uczak and Schoen~\cite{lucz-scho04}, Nathanson~\cite{nath03a}).

The study of representation functions for sets of nonnegative integers is more complicated.  In this case, every integer has only finitely many representations.  It is an open problem to describe the set of functions $f: \N_0 \rightarrow \N_0$ that are representation functions.  A special case of this inverse problem is one of the most famous problems in additive number theory: The conjecture of Erd\H os and Tur\' an that the representation function of an asymptotic basis for the nonnegative integers must be unbounded.

In the last few years there has been considerable work on inverse problems for representation functions (for example,~\cite{borw-choi-chu06, chen07, nath07f, nath07h, grek-hadd-helo-pihk03, hadd-helo04, lee07a, nath04e,  nest-serr04, tang-chen07}).

\section{Bases associated to linear forms}
Let $F(x_1,\ldots,x_m) = u_1 x_1 + \cdots + u_mx_m$ be an $m$-ary linear form with nonzero, relatively prime integer coefficients $u_1, \ldots, u_m$.    Let $A_1, \ldots, A_h$ be sets of integers.  We define the set 
\[
F(A_1,\ldots, A_m) = \{F(a_1,\ldots,a_m) : a_i \in A_i \text{ for $i=1,\ldots,m$}\}.
\]
The \emph{representation function} associated with the form $F$ is 
\[
R_{A_1,\ldots,A_m,F}(n) = \card\left( \{ (a_1,\ldots,a_m)\in A_1\times \cdots \times A_m:  F(a_1,\ldots, a_m) = n \} \right).
\]
This is a function from \Z\ into $\N_0 \cup \{\infty\}.$
If $A_i=\emptyset$ for some $i=1,\ldots,m,$ then $F(A_1,\ldots, A_m) = \emptyset$ and $R_{A_1,\ldots,A_m,F}(n) = 0$ for all $n \in \Z.$  

If $A_i = A$ for all $i=1,\ldots,m,$ then we write
\[
F(A) = F(A,\ldots,A) = \{F(a_1,\ldots, a_m) : a_i \in A \text{ for $i=1,\ldots,m$}\}
\]
and
\[
R_{A,F}(n) = R_{A,\ldots,A,F}(n) 
= \card\left( \{ (a_1,\ldots,a_m)\in A^m :  F(a_1,\ldots, a_m) = n \} \right).
\]
The sumset $mA$ is a set of the form $F(A),$ where $F$ is the linear form $F(x_1,\ldots,x_m) = x_1 + x_2 + \cdots + x_m$.

Let $S$ be a set of integers and let $F$ be an $m$-ary linear form.  The set $A$ is a \emph{basis for $S$ with respect to $F$} if $R_{A,F}(n) \geq 1$ for all integers $n \in S.$
The set $A$ is an \emph{basis with respect to $F$ for almost all $S$} if  $\{ n \in S : R_{A,F}(n) = 0\}$ has asymptotic density zero.

Many classical problems about bases in additive number theory have natural analogues for bases with respect to an $m$-ary linear form.  For example, a \emph{direct problem} for representation functions is:  Given a form $F$ and a set $A$ of integers, compute the representation function $R_{A,F}.$  An \emph{inverse problem} for representation functions is:  Given a form $F$ and a function $f:\Z \rightarrow \N_0 \cup \{\infty\},$ does there exist a set $A$ of integers such that $R_{A,F}(n) = f(n)$ for all $n\in \Z$?  
In particular, if $f(n)=1$ for all integers $n$, does there exist a set $A$ such that  $R_{A,F}(n) = 1$  for all integers $n$?  Such a set is called 
a \emph{unique representation basis with respect to the form $F$}.

In this paper we solve the inverse problem for bases with respect to a binary linear form.  Let $F(x_1,x_2) = u_1 x_1 + u_2x_2$ be a binary linear form whose coefficients are nonzero, relatively prime  integers.   We shall prove that if $u_1u_2 \neq \pm 1, -2$ and if $f:\Z \rightarrow \N_0 \cup\{\infty\}$ is any function such that $f^{-1}(0)$ has zero density, then there exists a set $A$ of integers such that $R_{A,F}=f.$  Equivalently, every nonzero function $f$ such that the set $f^{-1}(0)$ has asymptotic density zero  is the representation function of a basis for almost all  \Z\ with respect to $F$ for every binary linear form $F(x_1,x_2) \neq x_1 \pm x_2$  or $2x_1-x_2.$

Related work on the additive number theory of finite sets defined by linear forms appears in~\cite{nath07j,nath07d}.

\section{Construction of bases for binary linear forms}
We begin with two simple observations.

\bl     \label{bfbf:lemma:u}
Let $u_1$ and $u_2$ be nonzero, relatively prime integers with $u_1u_2 \neq \pm 1.$  The following four integers are pairwise distinct:
\[
-u_1^2, -u_1u_2, u_1 u_2, u_2^2
\]
and the following three integers are pairwise distinct:
\[
u_2^2-u_1^2, u_2(u_1+u_2),-u_1(u_1+u_2).
\]
If, in addition, $u_1u_2 \neq -2,$ then these seven integers are pairwise distinct. 
\el

\bpf
This is a straightforward verification.
\epf

\bl       \label{bfbf:lemma:density}
Let $A_1,\ldots, A_n$ be sets of integers and let $d_U(A_i)=\alpha_i$ for $i=1,\ldots, n.$  Then $d_U\left( \bigcup_{i=1}^n A_i \right) \leq \sum_{i=1}^n \alpha_i.$  In particular, the union of a finite number of sets of asymptotic density zero has asymptotic density zero.

If $A$ has asymptotic density zero and if $r,s,t$ are integers with $r\neq 0,$ then the set $\{n\in \Z : rn \in \{ sa+t : a\in A\} \}$ has asymptotic density zero.
\el

\bpf
For every $\varepsilon > 0$ there is a number $x_0(\varepsilon)$ such that  $A_i(-x,x) \leq (\alpha_i + \varepsilon/n)(2x+1)$ for all $x \geq x_0(\varepsilon)$.  If $A =  \bigcup_{i=1}^n A_i$,  then 
\[
A(-x,x) \leq \sum_{i=1}^nA_i(-x,x) \leq  \sum_{i=1}^n\left( \alpha_i + \varepsilon/n \right) (2x+1) = \left( \sum_{i=1}^n\alpha_i + \varepsilon\right) (2x+1) 
\]
for all $x \geq x_0(\varepsilon)$, and so $d_U(A) \leq \sum_{i=1}^n\alpha_i.$  In particular, if $\alpha_i = 0$ for all $i = 1,\ldots, n,$ then $d(A) = d_U(A) = 0.$  

Finally, subsets, dilations, and translations of sets with asymptotic density zero also have asymptotic density zero.  The last statement of the Lemma follows from this observation.
\epf

The following result is fundamental.

\bl  \label{bfbf:lemma:fundamental}
Let $F(x_1,x_2) = u_1 x_1+u_2 x_2$ be a binary linear form whose coefficients $u_1, u_2$ are nonzero, relatively prime integers with $u_1u_2 \neq \pm 1$ and $u_1u_2 \neq -2.$   Let $W$ be a set of integers with asymptotic density zero.  Let $A'$ be a finite set of integers and let $b$ be any integer such that the sets $W$ and $F(A') \cup \{ b \}$ are disjoint.   There exists a set $C$ with $A'\subseteq C $ and $|C \setminus A'| = 2$ such that
\beq  \label{bfbf:rArC}
R_{C,F}(n) = 
\begin{cases}
R_{A',F}(b) + 1 & \text{if $n = b$} \\
R_{A',F}(n)       & \text{if $n\in F(A')\setminus \{b \}$} \\
1         & \text{if $n\in F(C)\setminus \left( F(A') \cup \{ b\}  \right)$} \\
0         & \text{if $n \in W$.}  
\end{cases}
\eeq
\el

\bpf
Since $\gcd(u_1,u_2)=1,$ there exist integers $v_1$ and $v_2$ such that $F(v_1,v_2) = u_1v_1+u_2v_2 = 1.$  For every integer $t$ we have
\begin{align*}
F(bv_1+u_2t,bv_2-u_1t) & = u_1(bv_1+u_2t)+u_2(bv_2-u_1t) \\
& = b(u_1v_1+u_2v_2) = b.
\end{align*}
We introduce the sets 
\[
B_t =  \{ bv_1+u_2t,bv_2-u_1t \} 
\]
and 
\[
C_t = A' \cup B_t.
\]
Note that the conditions $\gcd(u_1,u_2) = 1$ and $u_1u_2 \neq 0,\pm 1$ imply that $u_1\pm u_2 \neq 0.$
If $(u_1+u_2)t\neq b(v_2-v_1),$ then $bv_1+u_2t \neq bv_2-u_1t$ and $|B_t|=2.$   Similarly, $A' \cap B_t \neq \emptyset$ if and only if 
$u_2t = a - bv_1$ or $u_1t = bv_2-a$ for some $a\in A'.$  Since the set $A'$ is finite, it follows that  $A' \cap B_t = \emptyset$ and 
$|C_t| = |A' \cup B_t| = |A'|+2$ for all but finitely many integers $t$.

We shall prove that there exist infinitely many integers $t$ such that  the set $C_t$ also satisfies conditions~\eqref{bfbf:rArC}.  Note that $F(C_t)$ is the union of four sets:
\[
F(C_t) = F(A')  \cup F(A',B_t)  \cup F(B_t,A')  \cup F(B_t).
\]
We have
\[
F(A',B_t) = \{ F(a,bv_1+u_2t) : a\in A'\} \cup \{ F(a, bv_2-u_1t ) : a \in A'\}.
\]
If $x \in  \{ F(a,bv_1+u_2t) : a\in A'\},$ then there exists $a\in A'$ such that
\[
x = F(a,bv_1+u_2t) = u_1 a + u_2(bv_1+u_2t) = (u_1a + u_2v_1b) + u_2^2t.
\]
If $x \in \{ F(a, bv_2-u_1t ) : a \in A'\},$ then there exists $a\in A'$ such that
\[
x =  F(a, bv_2-u_1t ) =  (u_1a + u_2v_2b) - u_1 u_2 t.
\]
For every integer $t$, the functions $F(a,bv_1+u_2t)$ and $F(a,bv_2-u_1t )$ are  strictly monotonic functions of $a$.
If  
\[
\{ F(a,bv_1+u_2t) : a\in A'\} \cap \{ F(a, bv_2-u_1t ) : a \in A'\} \neq \emptyset
\]
then there exist $a,a' \in A'$ such that 
\[
(u_1a + u_2v_1b) + u_2^2t = (u_1a' + u_2v_2b) - u_1 u_2 t
\]
or, equivalently,
\[
u_2(u_1+u_2) t = u_1(a' - a) + u_2( v_2 - v_1)b.
\]
Since $u_2(u_1+u_2) \neq 0$ and the set $A'$ is finite, it follows that 
for all but finitely many integers $t$ we have 
\[
\{ F(a,bv_1+u_2t) : a\in A'\} \cap \{ F(a, bv_2-u_1t ) : a \in A'\}  = \emptyset
\]
for all $a,a' \in A',$ and so
\[
R_{A',B_t,F}(n) \leq 1  \qquad\text{for all $n\in \Z.$}
\]

The set $F(A') \cup \{b\} \cup W$ has zero asymptotic density.  If 
\[
F(A', B_t) \cap \left( F(A') \cup \{b\} \cup W \right) \neq \emptyset
\]
then either 
\[
 (u_1a + u_2v_1b) + u_2^2t \in  F(A') \cup \{b\} \cup W
 \]
for some $a\in A'$ or 
\[
(u_1a' + u_2v_2b) - u_1 u_2 t \in  F(A') \cup \{b\} \cup W
\]
for some $a'\in A'$.  In both cases, by Lemma~\ref{bfbf:lemma:density}, the set of integers $t$ for which the membership relation is possible is a set of integers of asymptotic density zero.  Equivalently, except for a set  of integers $t$ of asymptotic density zero, we have 
\[
F(A', B_t) \cap \left( F(A') \cup \{b\} \cup W \right) = \emptyset
\]
and
\[
R_{A',B_t,F}(n) = 0  \qquad\text{for all $n\in F(A') \cup \{b\} \cup W .$}
\]

Similarly,
\begin{align*}
F(B_t, A') 
& = \{ F(bv_1+u_2t,a) : a\in A'\} \cup \{ F( bv_2-u_1t, a ) : a \in A'\} \\
& =  \{  (u_1v_1 b + u_2 a) + u_1 u_2 t  : a\in A'\} 
\cup \{ (u_1v_2b + u_2 a) - u_1^2 t  : a \in A'\}
\end{align*}
and, except for a set  of integers $t$ of asymptotic density zero,
\[
R_{B_t,A',F}(n) \leq 1  \qquad\text{for all $n\in \Z$}
\]
and
\[
R_{B_t, A', F}(n) = 0  \qquad\text{for all $n\in F(A') \cup \{b\} \cup W .$}
\]

If $F(A', B_t) \cap F(B_t, A') \neq \emptyset,$ then there exist integers $a,a' \in A'$  that satisfy at least one of the following four equations:
\begin{align*}
(u_1a + u_2v_1b) + u_2^2t & = (u_1v_1 b + u_2 a') + u_1 u_2 t \\
(u_1a + u_2v_1b) + u_2^2t & = (u_1v_2 b + u_2 a') - u_1^2 t  \\
(u_1a + u_2v_2b) - u_1 u_2 t & = (u_1v_1 b + u_2 a') + u_1 u_2 t \\
(u_1a + u_2v_2b) - u_1 u_2 t & = (u_1v_2 b + u_2 a') - u_1^2 t .
\end{align*}
Equivalently, 
\begin{align*}
u_2 ( u_2 - u_1) t & =  u_2 a' - u_1a  +  ( u_1- u_2) v_1b \\
(u_1^2  + u_2^2) t & =   u_2 a' - u_1a + (u_1v_2  -  u_2v_1) b \\
 -2 u_1 u_2 t & =  u_2 a' - u_1a + (u_1v_1  - u_2v_2)b \\
 u_1( u_1 - u_2) t & = u_2 a' - u_1a + (u_1 - u_2) v_2 b  .
\end{align*}
Since the coefficients of $t$ are nonzero and the set $A'$ is finite, it follows that there are only finitely many integers $t$ that can satisfy at least one of these equations for some $a,a' \in A$.  Except for this finite set of $t$, we have
\[
F(A', B_t) \cap F(B_t, A') = \emptyset
\]
and
\[
R_{A', B_t,  F}(n)  + R_{B_t, A', F}(n) \leq 1  \qquad\text{for all $n \in \Z.$}
\]

Finally, we consider the set $F(B_t),$ which consists of the integer $b$ and the three integers 
\[
(u_1v_2+u_2v_1)b+(u_2^2-u_1^2)t
\]
\[
(u_1+u_2)v_1b+u_2(u_1+u_2)t
\]
\[
(u_1+u_2)v_2b -u_1(u_1+u_2)t .
\]
The coefficients of $t$ in the last three expressions are distinct nonzero integers.    This implies that $|F(B_t)|=4$ for all but finitely many $t$, and also that
\[
\left( F(B_t) \setminus \{ b \} \right) \bigcap \left( F(A') \cup W \right)  = \emptyset 
\]
except for certain integers $t$ belonging to a set of asymptotic density zero.

The coefficients of $t$ in the integers in $F(A',B_t) \bigcup F(B_t,A')$ are $u_2^2$, $\pm u_1u_2, -u_1^2.$   The coefficients of $t$ in the integers in $F(B_t)\setminus \{ b\}$ are $u_2^2-u_1^2$, $u_2(u_1+u_2)$, $-u_1(u_1+u_2).$   
Since $u_1u_2\neq -2$,  these seven numbers are pairwise distinct by Lemma~\ref{bfbf:lemma:u}, and so the sets  $ F(A',B_t) \cup F(B_t,A')$ and $F(B_t)\setminus \{ b\}$  are pairwise disjoint for all but finitely integers $t$.  Since the union of a finite number of sets of asymptotic density zero is still a set of asymptotic density zero, it follows that, for almost all integers $t$, the set $C_t$ satisfies the requirements of the Lemma.  This completes the proof.
\epf

\bt\label{bfbf:theorem:inverse} 
Let $F(x_1,x_2) = u_1 x_1+u_2 x_2$ be a binary linear form whose coefficients $u_1, u_2$ are nonzero, relatively prime integers such that $u_1u_2 \neq \pm 1$ and $u_1u_2 \neq -2.$   Let $f: \Z \rightarrow \N_0 \cup \{\infty\}$ be any function such that the set $f^{-1}(0)$ has asymptotic density zero.  There exists a set $A$ of integers such that  $R_{A,F}(n) = f(n)$ for all integers $n$.
\et

\bpf
Let $W=f^{-1}(0)$ and let $\{b_i\}_{i=1}^{\infty}$ be a sequence of integers that 
\[
\left| \{ i\in \N : b_i = n \}\right| = f(n)  \qquad\text{for all $n \in \Z$}.
\]
In particular, $b_i \notin W$ for all $i\in \N$.  

Let $A_0 = \emptyset$.  Then $W \cap \left( F(A_0)\cup \{ b_1\}\right) = \emptyset$.
Applying Lemma~\ref{bfbf:lemma:fundamental} with $A'=A_0$ and $b=b_1,$  we obtain a set $C_t$ such that $R_{C_t,F}(n)=1$ if $n\in F(C_t),$ and $F(C_t)\cap W = \emptyset.$  Let $A_1 = C_t.$    
Then $R_{A_1, F}(n)  \leq f(n)$ for all $n\in \Z.$  

Let $i \geq 2$ and suppose that we have constructed sets $A_0 \subseteq A_1 \subseteq \cdots \subseteq A_{i-1}$ with $R_{A_{i-1}, F}(n)  \leq f(n)$ for all $n\in \Z.$    
If $R_{A_{i-1},F}(b_i)=  f(b_i)$, then let $A_i = A_{i-1}$.
Suppose that $R_{A_{i-1},F}(b_i) < f(b_i)$.
Since $W \cap \left( F(A_{i-1})\cup \{ b_i\} \right) = \emptyset$,
we can apply 
Lemma~\ref{bfbf:lemma:fundamental} 
with $A'=A_{i-1}$ and $b = b_i.$ We obtain a set $C_t$ satisfying conditions~\eqref{bfbf:rArC}.  Let $C_t = A_i$.  This procedure gives an infinite increasing sequence of finite sets $A_0 \subseteq A_1 \subseteq A_2 \subseteq \cdots$ 
such that, for all $n \in \Z,$ we have $R_{A_i,F}(n) \leq f(n)$ for all $i\in \N$ and $\lim_{i\rightarrow \infty} R_{A_i,F}(n)=  f(n)$.
It follows that the set $A = \bigcup_{i=0}^{\infty}A_i$ satisfies $R_{A,F}(n)= f(n)$ for all $n\in \Z.$   This completes the proof.
\epf

\bt\label{bfbf:theorem:URB} 
Let $F(x_1,x_2) = u_1 x_1+u_2 x_2$ be a binary linear form whose coefficients $u_1, u_2$ are nonzero, relatively prime integers such that $u_1u_2 \neq \pm 1$ and $u_1u_2 \neq -2.$  There exists a unique representation basis with respect to the form $F$.
\et

\bpf
Apply Theorem~\ref{bfbf:theorem:inverse}  to the function $f(n)=1$ for all $n \in \Z.$
\epf

\section{Sidon sets}
The set $A$ will be called a \emph{Sidon set} with respect to the $m$-ary linear form $F$ if $R_{A,F}(n) \leq 1$ for all integers $n$.  More generally, we shall call the set $A$ a $B_F[g]$-set with respect to the  form $F$ if $R_{A,F}(n) \leq g$ for all integers $n$. With respect to the classical additive form $F(x_1,\ldots, x_m) = x_1 + \cdots + x_m,$ Erd\H os and Tur\' an conjectured that no $B_F[g]$-set of nonnegative integers is an asymptotic basis for the nonnegative integers.    The analogue of this conjecture for arbitrary $m$-ary linear forms is not true.  Here is a simple example of a Sidon set of  nonnegative integers that is a basis for the nonnegative integers with respect to an $m$-ary form.

\bt
For $m \geq 2$ and $g \geq 2$, let 
\[
F(x_1,\ldots,x_m) = x_1+gx_2 + g^2 x_3 + \cdots + g^{m-1}x_m.
\]
Let $A$ be the set of all nonnegative integers whose $g^m$-adic representations 
use only the digits $\{0,1,\ldots,g-1\},$ that is,
\[
A = \left\{ \sum_{i=0}^{\infty} d_ig^{im} : d_i \in \{0,1,\ldots,g-1\} \text{ and $d_i = 0$ for all sufficiently large $i$}
\right\}.
\]
Then $A$ is a Sidon basis for $\N_0$ with respect to the form $F$.
\et

\begin{proof}
This follows immediately from the uniqueness of the $g$-adic representation of a  nonnegative integer.
\end{proof}

\def\cprime{$'$}
\providecommand{\bysame}{\leavevmode\hbox to3em{\hrulefill}\thinspace}
\providecommand{\MR}{\relax\ifhmode\unskip\space\fi MR }
\providecommand{\MRhref}[2]{%
  \href{http://www.ams.org/mathscinet-getitem?mr=#1}{#2}
}
\providecommand{\href}[2]{#2}

\end{document}